\documentclass[1p,draft]{my_elsarticle}
\usepackage[utf8]{inputenc}
\DeclareUnicodeCharacter{012D}{\u{\i}}
\DeclareUnicodeCharacter{012D}{\u\i}
\usepackage{textcomp} 
\usepackage{eurosym}  

\usepackage{amsmath,amssymb}
\usepackage{amstext,amsthm}
\usepackage{latexsym}
\usepackage{amsopn}
\usepackage{xspace}
\usepackage{stmaryrd} 

\usepackage{enumerate}
\usepackage[shortlabels]{enumitem}
\setlist[enumerate]{align=left,leftmargin=*,labelsep*=5pt}

\DeclareMathOperator{\spann}{span}

\DeclareMathOperator*{\slim}{s-lim}

\newcommand{\bgl}[1]{\begin{equation}\label{#1}}
\newcommand{\egl}{\end{equation}}

\newcommand{\Z}{\mathbb{Z}}

\newcommand{\eps}{\varepsilon}
\newcommand{\R}{\mathbb{R}}
\newcommand{\C}{\mathbb{C}}

\newcommand{\pbig}[1]{\ensuremath{\bigl( #1 \bigr)}}

\newcommand{\br}[1]{\ensuremath{\{ #1 \}}}

\newcommand{\brbigg}[1]{\ensuremath{\biggl\{ #1 \biggr\}}}

\newcommand{\abs}[1]{\ensuremath{\lvert #1 \rvert}}

\newcommand{\Oh}{\mathcal O}
\newcommand{\oh}{o}
\newcommand{\ds}{\displaystyle}

\newcommand{\du}{\,du}

\newcommand{\dS}{\,dS}

\newcommand{\Rae}{R^{\{\alpha\}}_{2j,\eps}\kern-1pt}

\newcommand{\CN}{\mathcal{N}}

\newcommand{\CD}{\mathcal{D}}

\newcounter{zaehler}

\def\neuerlabel{(\roman{zaehler})\hfil}
\newdimen\azlabelsep
\newdimen\aztopsep
\newdimen\azitemsep
\newenvironment{aufzaehlung}[1]%
{\begin{list}{}%
{\usecounter{zaehler}%
\settowidth{\labelwidth}{#1}
\leftmargin\labelwidth
\labelsep\azlabelsep
\addtolength{\leftmargin}{\labelsep}
\topsep\aztopsep
\itemsep\azitemsep
\renewcommand{\makelabel}[1]{\ifx##1\empty\neuerlabel\else ##1\fi}}}%
{\end{list}}

\newcommand{\ie}{i.\,e.\xspace}
\newcommand{\eg}{e.\,g.\xspace}

\theoremstyle{plain}
\newtheorem{theorem}{Theorem}[section]

\theoremstyle{definition}

\theoremstyle{remark}

\begin{document}

\title{\large The Third Boundary Value Problem of Potential Theory for the Exterior Ball and the Approximation behaviour of the solution;\\ a Novel Open Problem}

\author[plb]{P.\,L.~Butzer}
\ead{Butzer@rwth-aachen.de}

\author[plb]{R.\,L.~Stens\corref{cor1}}
\cortext[cor1]{Corresponding author}
\ead{stens@matha.rwth-aachen.de}

\address[plb]{Lehrstuhl A für Mathematik, RWTH Aachen, 52056 Aachen, Germany}

\begin{abstract}
The paper is concerned  with the interconnection of the boundary behaviour of the solutions of the exterior Dirichlet and Neumann problems of harmonic analysis for the unit ball  in $\R^3$ with the corresponding behaviour of the associated ergodic inverse problems for the punched unlimited space. The basis is the theory of semigroups of linear operators mapping a Banach space $X$ into itself. The rates of approximation play a basic role.

Another tool is a Drazin-like inverse operator $B$ for the infinitesimal generator $A$ of a semigroup that arises naturally in ergodic theory. This operator $B$ is a closed, not necessarily bounded, operator. It was introduced in a paper with U.~Westphal (1970/71) \cite{Butzer-Westphal_1970} and extended to a generalized setting with J.\,J.~Koliha (2009) \cite{Butzer-Koliha_2009}.

The novel open problem concerns the third or Robin's problem of potential theory, the solution of which is not a semigroup of operators. Hence, the semigroup methods applied to Dirichlet's or Neumann's problem cannot be applied. The authors give several hints how to overcome these difficulties.

\end{abstract}


\maketitle

\section{Dirichlet's problem for the three dimensional unit ball and its interconnections with the associated inverse problem in unlimited space}\label{sec_Dirichlet}

Dirichlet's problem, also known as the first problem of potential theory,
%
%
%
is to determine a function $w(\varphi,\theta,r)$ twice continuously differentiable on its domain $[0,\pi]\times\R \times (0,\infty)$ which is $2\pi$-periodic with respect to $\theta$, and satisfies Laplace's equation (in spherical coordinates)
\begin{equation}\label{eq_Laplace_3d_kugel}
   \triangle w = \frac{\partial^2 w}{\partial r^2} +\frac2r\frac{\partial w}{\partial r}+ \frac1{r^2}\frac{\partial^2 w}{\partial \varphi^2}
   +\frac{\cos \varphi}{r^2\sin\varphi}\frac{\partial w}{\partial\varphi}
  +\frac{1}{r^2\sin^2\varphi}\frac{\partial^2 w}{\partial \theta^2}=0
\end{equation}
together with the boundary conditions
\begin{align}
    &\lim_{r\to 1+} \|w(\varphi,\theta,r)-h(\varphi,\theta)\|_{L^2(S)}=0,\label{eq_Laplace_3d_kugel_boundary_2}\\[1.5ex]
    &\lim_{r\to \infty} w(r,\varphi,\theta)=0\quad \big(\varphi\in[0,\pi],\theta\in[-\pi,\pi]\big),\label{eq_Laplace_3d_kugel_boundary_3}
\end{align}
where $h(\varphi,\theta)$ is a given $f\in L^2(S)$, and $S$ denotes the unit sphere in $\R^3$.

For the spherical coordinates $\varphi, \theta, r$ we use the convention
\[
   x= r \sin \varphi \cos\theta,\quad y=r \sin\varphi \sin \theta,\quad z=r\cos \varphi,
\]
where $x,y,z$ are the rectangular coordinates in $\R^3$, and $-\pi < \theta\le \pi$ is the longitude, while  $0\le \varphi\le \pi$ is the latitude on the sphere of radius $r=\sqrt{x^2+y^2+z^2}$.

The solution of the above problem can be written as a Fourier expansion with respect to the so-called (complex) spherical harmonics (see \cite{Lense_1953}, \cite{Petrowski_1955}, \cite[§§~16, 31, 32]{Triebel_1972}, \cite{Wang-Li_2006}, \cite{Atkinson-Han_2012}, \cite{Berens-Butzer-Pawelke_1968})
\begin{equation*}
  Y_k^m(\varphi,\theta):=\sqrt{\frac{(2k+1)(k-\abs{m})!}{4\pi(k+\abs{m})!}}P^{|m|}_k(\cos \varphi)e^{im\theta}\quad(k=0,1,2,\dots;\,m=0,\pm1,\dots,\pm k),
\end{equation*}
where
\begin{equation*}
  P^m_k(t)= (1-t^2)^{m/2}\frac{d^m}{dt^m}P_k(t)=\frac{(-1)^k(1-t^2)^{m/2}}{2^k k!}\frac{d^{k+m}}{dt^{k+m}} (1-t^2)^k\quad (m=0,1,\dots,k),
\end{equation*}
and $P_k=P_k^0$ are the Legendre polynomials defined by
\[
  P_k(t):=\frac{(-1)^k}{2^k k!}\frac{d^k}{dt^k} (1-t^2)^k.
\]
The $P^m_k$ are known as associated Legrendre functions. They are polynomials for $m$ even, and polynomials multiplied by a factor $(1-t)^{1/2}$ for $m$ odd.

It is well known that the spherical harmonics form a complete orthonormal system in the Hilbert space $L^2(S)$ with respect to the scalar product and the norm
\[
   <g_1,g_2>:= \iint_S g_1 \overline{g_2} \,dS
   = \int_{-\pi}^\pi \int_0^\pi g_1(\varphi,\theta)\overline{g_2(\varphi,\theta)}\sin \varphi\,d\varphi\, d\theta,\
  \|g\|_{L^2(S)}:= \sqrt{<g,g>},
\]
$\dS=\sin\varphi\,d\varphi\, d\theta$ being the surface area element of the unit sphere. In particular, there holds
\[
      \langle Y_k^m,Y_j^l\rangle=\delta_{k,j}\delta_{m,l}\quad(k=0,1,2,\dots;\,m=0,\pm1,\pm2\dots,\pm k).
\]
It follows that every function $g\in L^2(S)$ can be expanded into a Fourier series with respect to the complex spherical harmonics,
\[
  g(\varphi,\theta)= \sum_{k=0}^\infty \sum_{m=-k}^k \widehat g(m,k)Y_k^m(\varphi,\theta),
\]
the series being convergent in $L^2(S)$-norm. The Fourier coefficients $\widehat g(m,k)$ are given by
\begin{align*}
    \widehat g(m,k):=\langle g, Y_k^m\rangle=
    \int_{-\pi}^\pi \int_0^\pi g(\varphi,\theta)\overline{Y_k^m(\varphi,\theta)}\sin \varphi\,d\varphi\, d\theta.
\end{align*}

Now  the unique solution of Dirichlet's boundary value problem \eqref{eq_Laplace_3d_kugel} \eqref{eq_Laplace_3d_kugel_boundary_2}, \eqref{eq_Laplace_3d_kugel_boundary_3} is given by
\begin{equation*}\label{eq_Laplace_3d_kugel_solutionx}
    w(\varphi,\theta,r)= \sum_{k=0}^\infty r^{-(k+1)}\sum_{m=-k}^k \widehat h(m,k) Y_k^m(\varphi,\theta)\quad \big(\varphi\in[0,\pi],\theta\in\R, r>1\big).
\end{equation*}

Setting now $r=e^t$, \ie,
\begin{equation*}\label{eq_semi_3d}
  w(\varphi,\theta,r)= w(\varphi,\theta ,e^t) = \sum_{k=0}^\infty e^{-(k+1)t}\sum_{m=-k}^k \widehat h(m,k)Y_k^m(\varphi,\theta),
\end{equation*}
it follows that $T_A(t)h(\varphi,\theta):=w(\varphi,\theta,e^t)=w_A(\varphi,\theta,e^t)$ defines a semigroup of class $(C_0)$ with infinitesimal generator $A$ given by

\begin{equation*}\label{eq_generator_3d}
  A h(\varphi,\theta)=-\sum_{k=0}^\infty (k+1) \sum_{m=-k}^k \widehat h(m,k)Y_k^m(\varphi,\theta).
\end{equation*}
Its domain $\CD(A)$ can be characterized in terms of the Fourier coefficients, namely,
\begin{equation*}\label{eq_generator_domain_3d}
   \CD(A)=\brbigg{g\in L^2(S); \sum_{k=0}^\infty (k+1)^2\sum_{m=-k}^k \abs{\widehat g(m,k)}^2<\infty},
\end{equation*}
whereas the nullspace $\CN(A)$ is trivial, \ie, $A$ is injective. The inverse $A^{-1}$ is defined on all of $L^2(S)$ and has the Fourier series representation
\begin{equation*}
  A^{-1}h(\varphi,\theta)=-\sum_{k=0}^\infty (k+1)^{-1} \sum_{m=-k}^k \widehat h(m,k)Y_k^m(\varphi,\theta).
\end{equation*}
Furthermore, the semigroup generated by $A^{-1}$ is given by
\[
  V_{A^{-1}}(t)h(\varphi,\theta) =w_{A^{-1}}(e^t,\varphi,\theta)= \sum_{k=0}^\infty e^{-t/(k+1)}\sum_{m=-k}^k \widehat h(m,k)Y_k^m(\varphi,\theta).
\]

We also need the resolvent $R(\lambda;A)$ of the operator $A$, which has the representation in terms of the associated semigroup $T_A(t)$,
\begin{equation*}\label{eq_resolvent_representation}
  R(\lambda;A)h(\varphi,\theta)=\int_0^\infty e^{-\lambda u}T_A(u)h(\varphi,\theta)\du \quad(\lambda\in\C,\Re \lambda>0),
\end{equation*}
and likewise for $R(\lambda;A^{-1})$. It follows easily from \eqref{eq_resolvent_representation} that
\[
  \slim_{\lambda\to\infty} r_A(\varphi,\theta;\lambda;h)
  :=\slim_{\lambda\to\infty} \lambda R(\lambda;A)h(\varphi,\theta)
  =\slim_{t\to 0+}T_A(t)h(\varphi,\theta)=h(\varphi,\theta).
\]
See \cite[Section~1.3]{Butzer-Berens_1967} in this respect.

In the following, we often write $w_A(r,\varphi,\theta;h)$ and $w_{A^{-1}}(r,\varphi,\theta;h)$ for  $w_A(r,\varphi,\theta)$ and $w_{A^{-1}}(r,\varphi,\theta)$, respectively, in order to indicate the dependence of the boundary value $h$.

Concerning the order of approximation of $w(\varphi,\theta,r;h)$ towards $h(\varphi,\theta)$ as well as the ergodic behaviour of $w_{A^{-1}}(\rho,\varphi,\theta;h)$ one has by a general theorem of semigroup theory (see \cite{Gessinger_1998,Gessinger_Diss_1997,Butzer-Stens_2918}),

\begin{theorem}\label{thm_Dirichlet}
 Let $w_A(r,\varphi,\theta;h):=V_A(t)h(\varphi,\theta)$, $w_{A^{-1}}(r,\varphi,\theta;h):=V_{A^{-1}}(t)h(\varphi,\theta)$ with $r=e^t$. \medskip

 \noindent a) For $\alpha \in (0,1]$, the following six assertions are equivalent:

 \begin{aufzaehlung}{(iii)}
\item $\ds \|w_A(\varphi,\theta,r;h)-h(\varphi,\theta)\|_{L^2(S)}
=\begin{cases}
 \oh(\log r)\\[.5ex]
 \Oh\big((\log r)^\alpha\big)
 \end{cases}\quad(r \to 1+)$,\\[2ex]
\item $\|r_A(\varphi,\theta;\lambda;h)-h(\varphi,\theta)\|_{L^2(S)}=\begin{cases}
 \oh(\lambda^{-1})\\[.5ex]
 \Oh(\lambda^{-\alpha})
 \end{cases}\quad (\lambda\to\infty)$,
\item $\ds\bigg\|\frac1{\log r} \int_1^{r} w_{A^{-1}}(\rho,\varphi,\theta;h)\,\frac{d\rho}{\rho}\bigg\|_{L^2(S)}
=\begin{cases}
 \oh\big((\log r)^{-1}\big)\\[.5ex]
 \Oh\big((\log r)^{-\alpha}\big)
 \end{cases}
 \quad (r \to \infty)
 $,\\[2ex]
 \item $\|r_{A^{-1}}(\varphi,\theta;\lambda;h)\|_{L^2(S)}=\begin{cases}
 \oh(\lambda)\\[.5ex]
 \Oh(\lambda^{-\alpha})
 \end{cases}\quad (\lambda\to 0+)$,

\item $K(t,h;L^2(S),\CD(A)) =\begin{cases}
 \oh(t)\\[.5ex]
 \Oh(t^{\alpha})
 \end{cases}\quad (t\to 0+)$,


 \item $\begin{cases}h= 0\quad a.\,e.,\\[.5ex]
            h\in \CD(A) \text{\quad if } \alpha =1. \end{cases}$
\end{aufzaehlung}

\noindent b) There exist elements $h_\alpha, h_\alpha^*, h_\alpha^{**}\in L^2(S)$, such that

\begin{aufzaehlung}{(iii)}

%
%
\item  $\ds\bigg\|\frac1{\log r} \int_1^{r} w_{A^{-1}}(\rho,\varphi,\theta;h_\alpha)\,\frac{d\rho}{\rho}\bigg\|_{L^2(S)}\ds 
\begin{cases}
 =\Oh\big((\log r)^{-\alpha}\big)\\[.5ex]\neq\oh\big((\log r)^{-\alpha}\big)
 \end{cases}\quad(r \to \infty)
 $,
\item  $\ds \big\|w_{A^{-1}}(\rho,\varphi,\theta;h_\alpha^*)- h_\alpha^*(\varphi,\theta)\big\|_{L^2(S)}
\begin{cases}
 =\Oh\big((\log \rho)^{\alpha}\big)\\[.5ex]\neq\oh\big((\log \rho)^\alpha\big)
 \end{cases}\quad(\rho \to 1+)
 $,
%
\item $\ds\big\|r_{A^{-1}}(\varphi,\theta;\lambda^{-1};h_\alpha^{**})\big\|_{L^2(S)}
    =\big\|r_A(\varphi,\theta;\lambda;h_\alpha^{**})-h_\alpha^{**}\big\|_{L^2(S)}
=\begin{cases}=\Oh(\lambda^{\alpha})\\[.5ex]
             \neq \oh(\lambda^{\alpha})
\end{cases}\\\rule{0pt}{0pt}\hfill
(\lambda\to 0+)
$.
\end{aufzaehlung}
\end{theorem}

\section{Neumann's boundary value problem}\label{sec_Neumann}

Neumann's problem, or second problem of potential theory, is to determine a function $w(\varphi,\theta,r)$ satisfying conditions \eqref{eq_Laplace_3d_kugel} and \eqref{eq_Laplace_3d_kugel_boundary_3} above, but with the boundary condition \eqref{eq_Laplace_3d_kugel_boundary_2} replaced by
\begin{equation*}\label{eq_Laplace_3d_kugel_boundary_2_second}
\lim_{r\to 1+}
    \bigg\|\frac{\partial}{\partial r}w(\varphi,\theta,r) -h(\varphi,\theta)\bigg\|_{L^2(S)}=0,
\end{equation*}

Its unique solution is given by (see, \eg, \cite[§§~16, 31, 32]{Triebel_1972}, \cite{Lense_1953}, and for an exact formulation in a two-dimensional setting \cite[p.~287]{Butzer-Nessel_1971})
\begin{equation}\label{eq_Laplace_3d_kugel_solution_Neumann}
    w(\varphi,\theta,r)= -\sum_{k=0}^\infty \frac{r^{-(k+1)}}{k+1}\sum_{m=-k}^k \widehat h(m,k) Y_k^m(\varphi,\theta)\quad \big(\varphi\in[0,\pi],\theta\in\R,r>1\big).
\end{equation}
In order to examine the order of approximation in \eqref{eq_Laplace_3d_kugel_boundary_2_second}, observe that
\begin{align*}
   \frac{\partial}{\partial r} w(\varphi,\theta,r)\Big|_{r=e^t}
   ={} & \sum_{k=0}^\infty r^{-(k+2)}\sum_{m=-k}^k \widehat h(m,k) Y_k^m(\varphi,\theta)\Big|_{r=e^t} \\[2ex]
   ={} & \sum_{k=0}^\infty e^{-(k+2)t}\sum_{m=-k}^k \widehat h(m,k) Y_k^m(\varphi,\theta)
         \quad \big(\varphi\in[0,\pi],\theta\in\R,r>0\big).
\end{align*}
The right-hand side defines a $C_0$-semigroup $S(t)$ on $L^2(S)$, and its generator $A$ is given by
\begin{equation*}\label{eq_generator_3d_neu}
  A h(\varphi,\theta)=-\sum_{k=0}^\infty (k+2) \sum_{m=-k}^k \widehat h(m,k)Y_k^m(\varphi,\theta)
\end{equation*}
with domain
\begin{equation*}\label{eq_generator_domain_3d_neu}
   \CD(A)=\brbigg{g\in L^2(S); \sum_{k=0}^\infty (k+2)^2\sum_{m=-k}^k \abs{\widehat g(m,k)}^2<\infty}.
\end{equation*}
The associated resolvent operator is given by
\[
 R(\lambda;A)h(\varphi,\theta)=\sum_{k=0}^\infty\frac{1}{\lambda+k+2}\sum_{m=-k}^k \widehat h(m,k) Y_k^m(\varphi,\theta)
         \quad \big(\varphi\in[0,\pi],\theta\in\R,t>0\big).
\]

The generator $A$ is injective and the inverse $A^{-1}$ is defined on all of $L^2(S)$ and has the Fourier series representation
\begin{equation*}
  A^{-1}h(\varphi,\theta)=-\sum_{k=0}^\infty \frac1{k+2} \sum_{m=-k}^k \widehat h(m,k)Y_k^m(\varphi,\theta).
\end{equation*}
Furthermore, the semigroup generated by $A^{-1}$ is given by
\begin{equation}\label{eq_semi_A_hoch_-1_3-dim}
  V_{A^{-1}}(t)h(\varphi,\theta) =w_{A^{-1}}(\varphi,\theta, e^t)= \sum_{k=0}^\infty e^{-t/(k+2)}\sum_{m=-k}^k \widehat h(m,k)Y_k^m(\varphi,\theta),
\end{equation}
having the resolvent operator
\[
R(\lambda;A^{-1})h(\varphi,\theta)=\sum_{k=0}^\infty\frac{k+2}{\lambda(k+2)+1}\sum_{m=-k}^k \widehat h(m,k) Y_k^m(\varphi,\theta)
         \quad \big(\varphi\in[0,\pi],\theta\in\R,t>0\big).
\]

The counterpart of Theorem~\ref{thm_Dirichlet} now reads.

\begin{theorem}\label{thm_Neumann}
 Let $w_A(\varphi,\theta,r;h):=w(\varphi,\theta,r)$ be defined by \eqref{eq_Laplace_3d_kugel_solution_Neumann}, and $w_{A^{-1}}(\varphi,\theta,r;h):=w_{A^{-1}}(\varphi,\theta, r)$ as in \eqref{eq_semi_A_hoch_-1_3-dim} with $r=\log t$. Further, let $r_A(\varphi,\theta;\lambda;h)= \lambda R(\lambda;A)h(\varphi,\theta)$, and $r_A^{-1}(\varphi,\theta;\lambda;h)=\lambda R(\lambda;A^{-1})h(\varphi,\theta)$. \medskip

 \noindent a) The following six assertions are equivalent:

 \begin{aufzaehlung}{(iii)}
\item $\ds \Big\|\frac\partial{\partial r}w_A(\varphi,\theta,r;h)-h(\varphi,\theta)\Big\|_{L^2(S)}\\[2ex]
\llap{$={}$}\bigg\|\sum_{k=0}^\infty r^{-(k+2)}\sum_{m=-k}^k \widehat h(m,k) Y_k^m(\varphi,\theta) -h(\varphi,\theta)\bigg\|_{L^2(S)}
=\begin{cases}
 \oh(\log r)\\[.5ex]
 \Oh\big((\log r)^\alpha\big)
 \end{cases}(r \to 1+)$,\\[2ex]
\item $\ds\|r_A(\varphi,\theta;\lambda;h)-h(\varphi,\theta)\|_{L^2(S)}\\[2ex]
\llap{$={}$}\bigg\|\sum_{k=0}^\infty\frac{\lambda}{\lambda+k+2}\sum_{m=-k}^k \widehat h(m,k) Y_k^m(\varphi,\theta) -h(\varphi,\theta)\bigg\|_{L^2(S)}
=\begin{cases}
 \oh(\lambda^{-1})\\[.5ex]
 \Oh(\lambda^{-\alpha})
 \end{cases}\quad (\lambda\to\infty)$,\\[2ex]
\item $\ds\bigg\|\frac1{\log r} \int_1^{r} w_{A^{-1}}(\varphi,\theta, \rho;h) \,\frac{d\rho}{\rho}\bigg\|_{L^2_{2\pi}}$\\[2ex]
%
%
  $\ds\llap{$={}$}\bigg\|\frac1{\log r} \int_1^{r}\sum_{k=0}^\infty \rho^{-\frac1{k+2}-1} \!\!
  \sum_{m=-k}^k \widehat h(m,k)Y_k^m(\varphi,\theta)\, d\rho \bigg\|_{L^2(S)}
=\begin{cases}
 \oh\big((\log r)^{-1}\big)\\[.5ex]
 \Oh\big((\log r)^{-\alpha}\big)
 \end{cases}\\
 \rule{0pt}{0pt}\hfill (r \to \infty)$,

 \item $\|w_{A^{-1}}(\varphi,\theta,r;h)\|_{L^2(S)}$\\[2ex]
  $\ds\llap{$={}$}\bigg\|\sum_{k=0}^\infty\frac{\lambda}{\lambda+k+2}\sum_{m=-k}^k \widehat h(m,k) Y_k^m(\varphi,\theta)\bigg\|_{L^2(S)}
 =\begin{cases}
 \oh(\lambda)\\[.5ex]
 \Oh(\lambda^{-\alpha})
 \end{cases}\quad (\lambda\to 0+)$,\\[2ex]
 %
%
\item $K(t,h;L^2(S),\CD(A)) =\begin{cases}
 \oh(t)\\[.5ex]
 \Oh(t^{\alpha})
 \end{cases}\quad (t\to 0+)$,


 \item $\begin{cases}f(x)= 0\quad a.\,e.,\\[.5ex]
            f\in \CD(A) \text{\quad if } \alpha =1. \end{cases}$
\end{aufzaehlung}

\noindent b) For any $\alpha\in (0,1]$ there exist elements $h_\alpha, h_\alpha^*, h_\alpha^{**}\in L^2(S)$, such that

\begin{aufzaehlung}{(iii)}

%
%
\item  $\ds\bigg\|\frac1{\log r} \int_1^{r} w_A(\varphi,\theta,\rho;h_\alpha)\,\frac{d\rho}{\rho}\bigg\|_{L^2(S)}
\begin{cases}
 =\Oh\big((\log r)^{-\alpha}\big)\\[.5ex]\neq\oh\big((\log r)^{-\alpha}\big)
 \end{cases}\quad(r \to \infty)
 $,
\item  $\ds \big\|w_{A^{-1}}(\varphi,\theta,\rho;h_\alpha^*) - h_\alpha^*(\varphi,\theta)\big\|_{L^2(S)}
\begin{cases}
 =\Oh\big((\log \rho)^{\alpha}\big)\\[.5ex]\neq\oh\big((\log \rho)^\alpha\big)
 \end{cases}\quad(\rho \to 1+)
 $,
%
\item $\ds\big\|r_{A^{-1}}(\varphi,\theta;\lambda^{-1};h_\alpha^{**})\big\|_{L^2(S)}
    =\big\|r_A(\varphi,\theta;\lambda)-h_\alpha^{**}(\varphi,\theta)\big\|_{L^2(S)}
=\begin{cases}=\Oh(\lambda^{\alpha})\\[.5ex]
             \neq \oh(\lambda^{\alpha})
\end{cases}\\
\rule{0pt}{0pt}\hfill(\lambda\to 0+)
$.
\end{aufzaehlung}

\end{theorem}

\section{The third or Robin's problem of potential theory}
This problem of potential theory, also named after V.\,G.~Robin%
\footnote{Victor Gustave Robin(1855--1897), who was professor at the Sorbonne in Paris, is especially known for the
Robin boundary conditions. He was awarded the Prix Franc{\oe}ur for 1893 and 1897, the Prix Poncelet for 1895; see \cite{Gustafson-Abe_1998a,Gustafson-Abe_1998b}.},
is to determine a function $w(\varphi,\theta,r)$ defined on $(0,\pi)\times\R \times (0,\infty)$ which is $2\pi$-periodic with respect to $\theta$, twice continuously differentiable on its domain, and satisfies Laplace's equation \eqref{eq_Laplace_3d_kugel}
together with the boundary conditions
\begin{align}
%
    &\lim_{r\to 1+}
    \bigg\|\alpha w(\varphi,\theta,r)+ \beta \frac{\partial}{\partial r}w(\varphi,\theta,r) -h(\varphi,\theta)\bigg\|_{L^2(S)}=0,\label{eq_Laplace_3d_kugel_boundary_2_third}\\[1.5ex]
    &\lim_{r\to \infty} w(\varphi,\theta, r)=0\quad \big(\varphi\in[0,\pi],\theta\in[-\pi,\pi]\big),\label{eq_Laplace_3d_kugel_boundary_3_third}
\end{align}
where $h(\varphi,\theta)$ is a given function in $L^2(S)$, $S$ being the unit ball in $\R^3$.

The general solution of equation \eqref{eq_Laplace_3d_kugel} satisfying the boundary conditions \eqref{eq_Laplace_3d_kugel_boundary_2_third} and \eqref{eq_Laplace_3d_kugel_boundary_3_third} is given by
\[
    w(\varphi,\theta,r)= \sum_{k=0}^\infty r^{-(k+1)}\sum_{m=-k}^k c_{m,k} Y_k^m(\varphi,\theta).
\]

For the coefficients $c_{m,k}$, one has by  \eqref{eq_Laplace_3d_kugel_boundary_2_third},
\[
    \lim_{r\to 1+}\Big(\alpha r^{-(k+1)}+\frac d{dr} \beta r^{-(k+1)}\Big) c_{m,k}
    = \big( \alpha-\beta (k+1)\big) c_{m,k}= \widehat h(m,k)
\]
This yields the unique solution
\begin{equation*}
  w(\varphi,\theta,r)= \sum_{k=0}^\infty \frac{r^{-(k+1)}}{\alpha-\beta (k+1)}\sum_{m=-k}^k \widehat h(m,k) Y_k^m(\varphi,\theta)\quad (\varphi\in[0,\pi], \theta\in\R, r>1),
\end{equation*}
provided $\alpha$ and $\beta$ are such that the denominator $\alpha-\beta (k+1)$ does not vanish for any $k$, which is, \eg, the case if they have different sign; see, \eg, \cite[pp.~181, 182]{Lense_1953}.

Since we are interested in the order of approximation in \eqref{eq_Laplace_3d_kugel_boundary_2_third}, we may try to proceed as in Sections~\ref{sec_Dirichlet} and \ref{sec_Neumann}, and consider

\setlength{\multlinegap}{.5cm}
\begin{multline}\label{eq_no_semigroup}
  \alpha w(\varphi,\theta,r)+ \beta \frac{\partial}{\partial r}w(\varphi,\theta,r)\Big|_{r=e^t}\\[2ex]
%
  =\sum_{k=0}^\infty \frac{\alpha e^{-(k+1)t}-\beta(k+1)e^{-(k+2)t}}{\alpha-\beta (k+1)}\sum_{m=-k}^k \widehat h(m,k) Y_k^m(\varphi,\theta)\qquad(t>0).
\end{multline}

Unfortunately, the family of operators, defined by \eqref{eq_no_semigroup} does not posses the semigroup property, apart from the cases $\beta=0$ (Dirichlet problem), or $\alpha=0$ (Neumann problem). This means that one cannot apply the general theory developed in \cite{Butzer-Stens_2918} in order to deduce results corresponding to those of Theorems~\ref{thm_Dirichlet} and \ref{thm_Neumann}.

A counterpart of Theorem~\ref{thm_Dirichlet}\,(i)$\Leftrightarrow$(v) or Theorem~\ref{thm_Neumann}\,(i)$\Leftrightarrow$(v), thus a direct and inverse approximation theorem for the operators \eqref{eq_no_semigroup}, may be proved by classical approximation theoretic methods including $K$-functional methods, where $\CD(A)$ has to be replaced by a suitable subspace of $L^2(S)$. Setting
\[
  U(r)h(\varphi,\theta):=  w(\varphi,\theta,r;h)
\]
the operators $\pbig{U(r)}_{r>1}$ mapping the Hilbert space $L^2(S)$ into $\spann\br{Y_k^m; k\in\Z, \abs m\le k}$ are clearly commutative, \ie, $U(r)U(s)=U(s)U(r)$. Hence the new approach presented in \cite{Butzer-Scherer_1972a, Butzer-Scherer_1972b} can be applied. For the Bernstein-type inequality needed for see, \eg, \cite{Mhaskar-Narcowich-Prestin-Ward_2010?,Dai-Xu_2013}. See also,  \cite[pp. ]{DeVore-Lorentz_1993}, \cite[pp.~]{Mhaskar_2000}. For the corresponding Jackson-type inequality see \cite{Butzer-Scherer_1972a,Butzer-Scherer_1972b,Dai-Ditzian_2008}, \cite[pp.~]{DeVore-Lorentz_1993}, \cite[pp.~]{Mhaskar_2000}, \cite[Chapter~5]{Wang-Li_2006}. For the basic spherical harmonic theory in question see, \eg, \cite{Wang-Li_2006,Atkinson-Han_2012,Dai-Xu_2013,Sansone_1959,Lense_1953}.

For the other equivalent assertions of Part~a), however, one has, first of all, to find a replacement for the infinitesimal generator $A$. With methods of spectral theory (see, \eg, \cite{Riesz-Nagy_1955,Edmunds-Evans_2018,Reed-Simon_1980}), which enables one to investigate the resolvent of $A$, one may be able to to deduce an equivalence of type (i)$\Leftrightarrow$(ii).

Concerning (iii) and (iv) of Theorems~\ref{thm_Dirichlet} or \ref{thm_Neumann}, the fundamental tool is the the so-called interconnection theorem, a link between the resolvents of $A$ and $A^{-1}$, where $A$ is the infinitesimal generator of a semigroup. It reads (see \cite{Butzer-Stens_2918}),
\begin{equation}\label{eq_interconnection}
   \lambda R(\lambda;A^{-1})f= f-\lambda^{-1} R(\lambda^{-1};A)f.
\end{equation}
The question is, whether there exists something similar, when $A$ is not the generator of a semigroup of operators?

The situation becomes even more complicated, when $A$ is not injective. In this case $A^{-1}$ has to be replaced by the generalized Drazin inverse $A^{\textrm{ad}}$, introduced by Butzer, Westphal and Koliha (see \cite{Butzer-Westphal_1970, Butzer-Westphal_1972,Butzer-Koliha_2009}). Here $A$ is a closed not necessarily bounded operator. If $A$ is the generator of a semigroup, then \eqref{eq_interconnection} holds for $A^{\textrm{ad}}$ instead of $A^{-1}$ in a slightly modified form; see \cite{Gessinger_1998}, \cite[p.~37]{Gessinger_Diss_1997},

Alternative approaches to counterparts of (iii) and (iv) of Theorems~\ref{thm_Dirichlet} or \ref{thm_Neumann} are using Shaw's theory of \textbf{A}-ergodic nets \cite{Shaw_1998} or via Hille's pseudoresolvents (see \cite[p.~215~ff.]{Yosida_1978}, \cite[p.~521~ff.]{Hille-Philipps_1957}).

\section*{Acknowledgement}
The authors wish to extend their thanks to Isaac Pesensen, Temple University, for his continued interest in our results, especially in connection with our paper \cite{Butzer-Stens_2918}. Whenever we sent him a new section he always gave us useful hints, in particular, with respect to the open problem.

\section{References}

\bibliographystyle{elsarticle-num}
\bibliography{literatur}

\end{document}